\documentclass{article}
\usepackage{amsmath}
\usepackage{amssymb}
\usepackage{amsthm}
\usepackage[toc,page]{appendix}
\usepackage{graphicx}
\numberwithin{equation}{section}
\usepackage{xcolor}

\newtheorem{theorem}{Theorem}

\title{On the Use of the Mellin Transform to Generate Families of Power, Hyperpower, Lambert and Dirichlet Type Series and Some Consequences}
\author{Larry Glasser{\footnote{Department of Physics, Clarkson University, Potsdam, NY, USA. laryg@clarkson.edu}} and Michael Milgram{\footnote{Consulting Physicist, Geometrics Unlimited Ltd., Deep River, Ontario Canada. mike@geometrics-unlimited.com}}}
\begin{document}\vskip 1in
\date{September 28, 2023}
\maketitle

Revised Sept. 17, 2024: corrected Eq. (3.41)
\vskip .2 in

\vskip 0.2in
\centerline{\bf ABSTRACT}\vskip .2in

This note is concerned with series of the forms $\sum f(a^n)$ and $\sum f(n^{-a})$ where f(a)
possesses a Mellin transform and $a > 1$ or $a<0$ respectively. Integral representations
are derived and used to transform these series in several ways yielding a selection of interesting integral evaluations involving Riemann's function $\zeta(s)$, limits and series representations containing hyperpowers. A number of examples of such sums are provided, each of which is investigated for possible new structure. In one case, we obtain a generalization of Riemann's classic relationship among the Zeta, Gamma and Jacobi Theta functions.

\vskip .2in

\section{Introduction} \label{sec:Intro}

Although series of the form
\begin{equation}
S=\sum_{n=0}^{\infty} f(a^n)
\label{eqno(1)}
\end{equation}

have been extensively studied as $q-$extensions of classical functions, specific examples only appear sporadically in mathematical tables. Some examples appear in the extensive tables of Prudnikov et.al. \cite[sections 5.4.11 -16]{prudnikov1}, as well as  Hansen \cite{Hansen} sections 11-13 and 17.9. General relationships between representative sums of the form studied here and similar products can also be found in \cite[Section 17]{Erdelyi3}.

In this note we assume that $f(x)$ possesses a Mellin transform
\begin{equation}
F(s) =\int_0^{\infty} x^{s-1} f(x) dx, \hspace{20pt} \Re(s)>0
\label{eqno(2)}
\end{equation}
and employ this property to obtain interesting identities, one of which can be connected to $q$-identities. Throughout, $j\,,n,k ,N$ are non-negative integers, $b\in\Re$ and $b>0$, all other variables are complex unless specified otherwise, $\gamma$ is the Euler-Mascheroni constant and $\psi_{q}(x)$ is the $q$-extension of the digamma function. The symbol $:=$ refers to symbolic replacement and we represent a frequently appearing sum in terms of a Jacobi-type analogue theta function
\begin{equation}
\omega(b,s)\equiv\overset{\infty}{\underset{j =1}{\sum}}\; {\mathrm e}^{-b\,j^{{s}}}\,.
\label{OmegId}
\end{equation}


Basically, we utilize Lebesgue's dominated convergence theorem 
\begin{equation}
S=\sum_{n=0}^{\infty}\int_{c-i\infty}^{c+i\infty}
\frac{ds}{2\pi i} a^{-ns} F(s)=\int_{c-i\infty}^{c+i\infty}\frac{ds}{2\pi i}\sum_{b=0}^{\infty} a^{-ns} F(s)
\label{eqno(3)}
\end{equation}
for some $c>0$ and $a>1$, ensuring the convergence of the series. Therefore, we have

\begin{theorem}
If $f$ possesses a Mellin transform $F$ and $Re\,( a)\,>1$, then
\begin{equation}
\sum_{n=0}^{\infty}f(a^n)=\frac{1}{2\pi i}\int_{c-i\infty}^{c+i\infty}\frac{F(s)}{1-a^{-s}}{ds}.
\label{eqno(4)}
\end{equation}
\end{theorem}
\noindent

Pursuing this principle of treating the free Mellin transform variable $a$ as the independent variable in a power series, we also consider its utilization in the form of a Dirichlet series; specifically, if

\begin{equation}
f \! \left(a \right) = 
\frac{1}{2\,\pi i} \int_{c -i\infty }^{c +i\infty }a^{-v}\,F \! \left(v \right)d v 
\label{Fbase}
\end{equation}
where $f(a)$ and $F(v)$ are related by an inverse Mellin transform, then, 

\begin{theorem} \label{sec:T2}
If $\Re(s)<0$ and the sum converges, we have
\begin{equation}
\overset{\infty}{\underset{j =1}{\sum}}\; f \! \left(j^{-s}\right)
 = 
\,\frac{1 }{2\,\pi i}\int_{c -i \infty }^{c +i\infty }\zeta \! \left(-s\,v \right) F \! \left(v \right)d v\,. 
\label{F0}
\end{equation}

\end{theorem}

{\bf Remark:} In another context \cite{MDet}, instead of summing or exponentiating, the variable $a$ was treated as a complex variable $a\equiv r\,\exp{i\theta}$, and the Mellin transforms were studied as a function of $\theta$.\newline

In the following Section \ref{sec:PowerS}, we present some examples based on Theorem 1. In the subsequent Section \ref{sec:Dirichlet}, we present some examples based on Theorem 2. From the first of these Sections, a few of the more interesting identities are listed below:

\begin{itemize}
\item{see \eqref{X1Minus}}
\begin{align*}
\overset{\infty}{\underset{n =0}{\sum}}\; \frac{{\mathrm e}^{-a^{n}}}{\sinh \! \left(a^{n}\right)}
 = &
\,-\frac{\ln \! \left(\pi \right)}{\ln \! \left(a \right)}-\frac{1}{2}+\frac{a}{a -1}+\frac{\gamma}{\ln \! \left(a \right)}+\overset{\infty}{\underset{k =1}{\sum}}\; \left(a^{k}-\coth \! \left(a^{-k}\right)\right) \\
&+\frac{4}{\ln \! \left(a \right)}\,\Re \! \left(\overset{\infty}{\underset{k =1}{\sum}}\; \zeta \! \left(\frac{2\,i\,\pi \,k}{\ln \left(a \right)}\right) \Gamma \! \left(\frac{2\,i\,\pi \,k}{\ln \left(a \right)}\right) 2^{-\frac{2\,i\,\pi \,k}{\ln \left(a \right)}}\right)\,;
\end{align*}

\item{see \eqref{Eq4}}
\begin{equation} \nonumber
\overset{\infty}{\underset{k =1}{\sum}}\; \frac{1}{1-a^{2\,k +1}}
 = \frac{1}{a -1}
-\frac{\psi_{{1}/{a^{2}}}\! \left(1\right)}{2\,\ln \! \left(a \right)}+\frac{\psi_{{1}/{a}}\! \left(1\right)}{\ln \! \left(a \right)}+\frac{\ln \! \left(\frac{a -1}{a +1}\right)}{2\,\ln \! \left(a \right)}\,;
\end{equation}

\item{see \eqref{E4B}}
\begin{equation} \nonumber
\overset{\infty}{\underset{j =1}{\sum}}\; \Re \! \left(\Gamma \! \left(\frac{2\,i\,\pi \,j}{\ln \! \left(b \right)}\right)\right)
 = 
\frac{\ln \! \left(b \right)}{2} \,\overset{\infty}{\underset{j =0}{\sum}}\! \left({\mathrm e}^{-{1}/{b^{j}}}-1+{\mathrm e}^{-b^{j}}\right)+\left(\frac{1}{4}-\frac{1}{2\,e}\right) \ln \! \left(b \right)+\frac{\gamma}{2}\,.
\end{equation}

From the next Section \ref{sec:Dirichlet}, interesting identities include

\item{see \eqref{F5b4}}
\begin{align*}
\int_{0}^{\infty}&\frac{{| \zeta \! \left(\frac{1}{2}+i\,v \right)|}^{2} }{\cosh \! \left(\pi \,v \right)} \left(\cos \! \left(v\,\ln \! \left(2\,\pi \right)\right) \cosh \! \left(\frac{\pi \,v}{2}\right)-\sin \! \left(v\,\ln \! \left(2\,\pi \right)\right) \sinh \! \left(\frac{\pi \,v}{2}\right)
\right) d v \\
 &= 
\sqrt{\pi}\left(\overset{\infty}{\underset{j =1}{\sum}}\; \frac{1}{{\mathrm e}^{j}-1}-\gamma\right);
\end{align*}

\item{see \eqref{Jb}}
\begin{equation} \nonumber
\zeta \! \left(-b\,s \right) \Gamma \! \left(b \right) = 
\int_{0}^{\infty}x^{b -1}\overset{\infty}{\underset{j =1}{\sum}}\; {\mathrm e}^{-x\,j^{-s}}d x,~~~s<0,
\end{equation}
reducing to Riemann's classic identities \cite[Eqs. (2.4.1) and (2.6.2)]{Titch2} when $s=-1$ and $s=-2$ respectively;

\item{see \eqref{T2bIdX}}
\begin{equation} \nonumber
\underset{n \rightarrow \infty}{\mathrm{lim}}\! \left(\overset{\infty}{\underset{j =1}{\sum}}\; \frac{1}{{\mathrm e}^{j^{{1}/{n}}}-1}-\zeta \! \left(n \right) \Gamma \! \left(n +1\right)\right)
 = \frac{1}{2-2\,{\mathrm e}} \,.
\end{equation}
\end{itemize}

\section{Examples based on Theorem 1} \label{sec:PowerS}
\subsection{Example \ref{sec:T1P1}} \label{sec:T1P1}

From \cite[Eq. 6.6(2)]{IxForms} with $f(x)=1/\sinh(ax), ~ F(s)=2\,a^{-s} \left(1-2^{-s}\right) \Gamma \! \left(s \right) \zeta \! \left(s \right)$,

we have the Mellin transform
\begin{equation}
\int_{0}^{\infty}\frac{x^{s -1}}{\sinh \! \left(a\,x \right)}d x = 
2\,a^{-s} \left(1-2^{-s}\right) \Gamma \! \left(s \right) \zeta \! \left(s \right)
\label{Erd2}
\end{equation}
yielding, after differentiating with respect to the variable $a$, 
\begin{equation}
\int_{0}^{\infty}\frac{x^{s}\,\cosh \! \left(a\,x \right)}{\cosh \! \left(2\,a\,x \right)-1}d x
 = 
a^{-s -1} \left(1-2^{-s}\right) \Gamma \! \left(s +1\right) \zeta \! \left(s \right),\quad Re\; s>1.\label{eqno(5)}
\end{equation}

According to \eqref{eqno(4)}, after summation, we find the inverse Mellin transform with $c>1$ 

\begin{equation}
\frac{1}{2\,\pi\,i}\int_{c -i\,\infty}^{c +i\,\infty }\frac{b^{-s -1} \left(1-2^{-s}\right) \Gamma \left(s +1\right) \zeta \left(s \right)}{1-a^{-s}}d s
 = 
\overset{\infty}{\underset{j =0}{\sum}}\; \frac{a^{j}\,\cosh \! \left(a^{j}\,b\right)}{\cosh \! \left(2\,a^{j}\,b\right)-1}
\label{Eq5b}
\end{equation}
giving the (possibly new) series
\begin{equation}
\sum_{j=0}^{\infty} 2^j \frac{\cosh(2^j\,b)}{\sinh^2(2^j\, b)}=\frac{2}{b^{2}}+2 \overset{\infty}{\underset{k =1}{\sum}}\; \frac{b^{2\,k -2}\,\zeta \! \left(1-2\,k \right)}{\Gamma \! \left(2\,k -1\right)}=\frac{1}{2}\,\frac{1}{\sinh^2(b/2)},
\label{eqno(6)}
\end{equation}
by evaluating the residues (see Appendix \ref{sec:Proof1}) as the contour is moved to negative infinity when $a=2$.

Cases like this, where the quantity $1-a^{-s}$ cancels from the denominator of (4) are scarce. However since the factor $1-a^{-s}=1-e^{-s\ln a} $ vanishes for complex $s$, when the contour is moved into the left-half $s-$plane, further interesting residue sums over the imaginary poles
$ 2\pi k i/\ln a$,     $k=0,\pm 1,\pm 2\cdots$
can be found, especially when $F(s)$ is meromorphic as will now be demonstrated.

\subsection{Example \ref{sec:T1p2}} \label{sec:T1p2}
Consider next \cite[Eq. 6.3(7)]{IxForms}
\begin{equation}
f(x)=\frac{1}{e^{x}-1},\quad F(s)=\Gamma(s)\zeta(s)\,\,, \hspace{10pt} \Re(s)>1.
\label{eqno(7)}
\end{equation} 

Thus,
\begin{equation}
\sum_{n=0}^{\infty}\frac{1}{e^{a^n}-1}=\frac{1}{2\pi i}\int_{c-i\infty}^{c+i\infty}\frac{\Gamma(s)\zeta(s)}{1-a^{-s}}ds\,, \quad c>1\,.
\label{eqno(8)}
\end{equation}

The integrand has a double pole at $s=0$, simple poles at  s=\,-k$, $ $k$ odd and $s=2k\pi i/\ln a$, $k\in{\cal N}$. Closing the contour to the left and summing the appropriate residues yields
\begin{equation}
\sum_{n=0}^{\infty}\frac{1}{e^{a^n}-1}=  
\frac{\gamma-\ln(2\pi \sqrt{a})}{\ln(a^2)}+\frac{a}{a-1}$$
$$+\sum_{k=0}^{\infty}\frac{\zeta(-2k-1)}{(2k+1)!(a^{2k+1}-1)} +\frac{2}{\ln a} {\cal R}\sum_{k=1}^{\infty}\Gamma\left(\frac{2\pi i k}{\ln a}\right)\zeta\left(\frac{2\pi i k}{\ln a}\right).
\label{eqno(9)}
\end{equation}

However, since
\begin{equation}
\sum_{k=0}^{\infty}\frac{\zeta(-2k-1)}{(2k+1)!(a^{2k+1}-1)} =-\frac{1}{2}\sum_{k=1}^{\infty}\,[\coth\left(a^{-k}/2\right)-2a^{k}]\,,\label{eqno(10)}
\end{equation}

(see Appendix \ref{sec:Proof2}), then
\begin{align} \nonumber
 {\cal R}\sum_{k=1}^{\infty}\Gamma\left(\frac{2\pi i k}{\ln a}\right)&\zeta\left(\frac{2\pi i k}{\ln a}\right)=
-\frac{a\ln a}{2(a-1)}-\frac{\gamma-\ln(2\pi \sqrt{a})}{4}\\
&+\frac{\ln a}{2}\sum_{k=0}^{\infty}\frac{1}{e^{a^k}-1} +
\frac{\ln a}{4}\sum_{k=1}^{\infty}\,\left(\coth\left(a^{-k}/2\right)-2a^{k}\right)\,.
\label{eqno(11)}
\end{align}

Similarly \cite[Eq. 6.3(6)]{IxForms}, which, by analogy to \eqref{eqno(7)}, provides the Mellin transform pair

\begin{equation}
f(x)=\frac{1}{e^{x}+1},\quad F(s)=(1-2^{(1-s)})\,\Gamma(s)\,\zeta(s)\,, \hspace{10pt} \Re(s)>0
\label{eqno(7a)}
\end{equation} 
leads to

\begin{align} \nonumber
\overset{\infty}{\underset{n =0}{\sum}}\; \frac{1}{{\mathrm e}^{a^{n}}+1}
 =& 
\,-\frac{1}{2\,\ln \! \left(a \right)}\ln \! \left(\frac{2}{\pi \,\sqrt{a}}\right)-\frac{1}{2}\,\overset{\infty}{\underset{k =1}{\sum}}\; \left(\coth \! \left({a^{-k}}/{2}\right)-2\,\coth \! \left(a^{-k}\right)\right)\\
&-\frac{\gamma}{2\,\ln \! \left(a \right)}-\frac{4}{\ln \! \left(a \right)}\,\Re \, \overset{\infty}{\underset{k =1}{\sum}}\; \left(2^{-\frac{2\,i\,\pi \,k}{\ln \left(a \right)}}-\frac{1}{2}\right) \Gamma \! \left(\frac{2\,i\,\pi \,k}{\ln \left(a \right)}\right) \zeta \! \left(\frac{2\,i\,\pi \,k}{\ln \left(a \right)}\right)\,.
\label{X2b}
\end{align}

By adding \eqref{X2b} and \eqref{eqno(9)} we obtain
\begin{align} \nonumber
\overset{\infty}{\underset{n =0}{\sum}}\; \frac{1}{\,\sinh \! \left(a^{n}\right)}&
 = 
\,-\overset{\infty}{\underset{k =1}{\sum}}\! \left(\coth \! \left({a^{-k}}/{2}\right)-\coth \! \left(a^{-k}\right)-a^{k}\right)-\frac{a}{1-a }-\frac{\ln \! \left(2\right)}{\ln \! \left(a \right)}\\
&+\frac{4}{\ln \! \left(a \right)}\,\Re \; \overset{\infty}{\underset{k =1}{\sum}}\! \left(1-2^{-\frac{2\,i\,\pi \,k}{\ln \left(a \right)}}\right) \Gamma \! \left(\frac{2\,i\,\pi \,k}{\ln \left(a \right)}\right) \zeta \! \left(\frac{2\,i\,\pi \,k}{\ln \left(a \right)}\right)\,,
\label{X1Add}
\end{align}

and by subtracting we find

\begin{align} \nonumber
\overset{\infty}{\underset{n =0}{\sum}}\; \frac{{\mathrm e}^{-a^{n}}}{\sinh \! \left(a^{n}\right)}
 = &
\,-\frac{\ln \! \left(\pi \right)}{\ln \! \left(a \right)}-\frac{1}{2}+\frac{a}{a -1}+\frac{\gamma}{\ln \! \left(a \right)}+\overset{\infty}{\underset{k =1}{\sum}}\; \left(a^{k}-\coth \! \left(a^{-k}\right)\right) \\
&+\frac{4}{\ln \! \left(a \right)}\,\Re \; \overset{\infty}{\underset{k =1}{\sum}}\; \zeta \! \left(\frac{2\,i\,\pi \,k}{\ln \left(a \right)}\right) \Gamma \! \left(\frac{2\,i\,\pi \,k}{\ln \left(a \right)}\right) 2^{-\frac{2\,i\,\pi \,k}{\ln \left(a \right)}}\,.
\label{X1Minus}
\end{align}

Note that in taking the limit of \eqref{X1Add} as $a\rightarrow 2$, the first series on the right-hand side is telescoping while the second series vanishes since $1-2^{-{2\,i\,\pi \,k}/{\ln \left(2\right)}} = 0$,  yielding
\begin{equation}
\overset{\infty}{\underset{n =0}{\sum}}\frac{1}{\sinh \! \left(2^{n}\right)}
 = \coth \! \left(\frac{1}{2}\right)-1\,,
\label{Ts2}
\end{equation}
an identity that could also be obtained by evaluating the listed identity \cite[Eq. 25.1.1]{Hansen}
\begin{equation}
\overset{N}{\underset{n =0}{\sum}}\;\frac{1}{\sin \! \left(2^{n}\,x\right)}
 = \cot \! \left(\frac{x}{2}\right)-\cot \! \left(2^{N}\,x \right)
\label{H25b}
\end{equation}
after setting $x=i$ with $N\rightarrow\infty$. See also \cite[Eq. 1.121.2]{G&R} and \cite{AMM11853}. It is notable that the left-hand side of \eqref{eqno(11)} is expressible in elementary terms when $a=2$.

\subsection{Example \ref{sec:T1p3}} \label{sec:T1p3}
In a more elementary vein, let us take 
\begin{equation}
f(a)=(a^2+1)^{-1},~~F(s)=(\pi/2)/\sin(\pi s/2), a>1, \quad 0<\Re(s)<2,
\end{equation}
yielding the identity

\begin{equation}
\overset{\infty}{\underset{k =0}{\sum}}\; \frac{1}{a^{2k}+1}
 = \,
-\frac{i}{4} \int_{1 -i\,\infty}^{1 +i\,\infty }\frac{1}{\sin \left({\pi \,s}/{2}\right) \left(1-a^{-s}\right)}d s\,. 
\label{X3}
\end{equation}

Except for the case $s=0$, the residues from the zeroes of the denominator term $1/(1-a^{-s})$ in the integrand of \eqref{X3} cancel, and closing the contour by transiting the poles $s=-2k, \quad k=1,2,3\cdots$ and the double pole $s=0$, we obtain the known \cite[Eq. (2.1a)]{schmidt2020catalog} Lambert series identity (originally attributed to Ramanujan)

\begin{equation}
\overset{\infty}{\underset{k =1}{\sum}}\; \frac{1}{a^{2\,k}+1} = 
\overset{\infty}{\underset{k =1}{\sum}}\; \frac{\left(-1\right)^{1+k}}{a^{2\,k}-1},\quad a>1
\label{X3a}
\end{equation}

which, by letting $a=\exp(x)$, is equivalent to
\begin{equation}
\sum_{k=1}^{\infty}(-1)^{k+1}\left({1}/{\tanh(kx)}-1\right)=\sum_{k=1}^{\infty}\left(1-\tanh(kx)\right) \quad (x>0)\label{eqno(15a)}
\end{equation}
because
\begin{equation}
\overset{\infty}{\underset{k =1}{\sum}}\;\left(\frac{{\mathrm e}^{k\,x}}{2\,\cosh \! \left(k\,x \right)}-1\right)
 = 
\overset{\infty}{\underset{k =1}{\sum}}\left(-1\right)^{1+k} \left(1-\frac{{\mathrm e}^{k\,x}}{2\,\sinh \! \left(k\,x \right)}\right)
\label{So}
\end{equation}
using \eqref{X3a}. The identity \eqref{X3a} provides an interesting connection to the $q-$digamma function, by considering the odd and even terms of each sum independently:

\begin{equation}
\overset{\infty}{\underset{k =1}{\sum}}\; \frac{1}{1-a^{2\,k +1}} = 
\,-\overset{\infty}{\underset{k =1}{\sum}}\; \frac{1}{1+a^{k}}+\overset{\infty}{\underset{k =1}{\sum}}\; \frac{1}{1-a^{2\,k}}-\frac{1}{1-a}
\label{Eq1}
\end{equation}
and
\begin{equation}
\overset{\infty}{\underset{k =1}{\sum}}\; \frac{1}{1-a^{2\,k +1}} = 
\overset{\infty}{\underset{k =1}{\sum}}\; \frac{1}{1-a^{k}}-\overset{\infty}{\underset{k =1}{\sum}}\; \frac{1}{1-a^{2\,k}}-\frac{1}{1-a}\,.
\label{Eq2}
\end{equation}
So, by equating the right-hand sides of each, we obtain
\begin{equation}
\overset{\infty}{\underset{k =1}{\sum}}\; \frac{1}{1+a^{k}} = 
\frac{1}{\ln \! \left(a \right)} \left(\ln \! \left({1 +1/a}\right)-\psi_{{1}/{a}}\! \left(1\right)+\psi_{{1}/{a^{2}}}\! \left(1\right)\right)
\label{Eq3}
\end{equation}
employing the identity \cite[Eq. (4)]{q-polyGamma}
\begin{equation}
\overset{\infty}{\underset{k =1}{\sum}}\; \frac{1}{1-a^{k}} = 
\frac{\psi_{{1}/{a}}\! \left(1\right)+\ln \! \left(a -1\right)}{\ln \! \left(a \right)}-1,\hspace{20pt}a>1,
\label{Pid1}
\end{equation}
where the $q-$digamma function is defined by \cite[Eq. (2)]{q-polyGamma}

\begin{equation}
\psi_{q}\! \left(z \right) = 
-\ln \! \left(1-q \right)+\ln \! \left(q \right) \overset{\infty}{\underset{k =0}{\sum}}\; \frac{q^{k +z}}{1-q^{k +z}}\,,\hspace{20pt} |q|<1\,.
\label{PsiDef}
\end{equation}
Either of \eqref{Eq1} or \eqref{Eq2} further yields
\begin{equation}
\overset{\infty}{\underset{k =1}{\sum}}\; \frac{1}{1-a^{2\,k +1}}
 = \frac{1}{a -1}
-\frac{\psi_{{1}/{a^{2}}}\! \left(1\right)}{2\,\ln \! \left(a \right)}+\frac{\psi_{{1}/{a}}\! \left(1\right)}{\ln \! \left(a \right)}+\frac{\ln \! \left(\frac{a -1}{a +1}\right)}{2\,\ln \! \left(a \right)}\,,
\label{Eq4}
\end{equation}
the alternating version of which is the well-known identity \cite[Eq. (5.1)]{schmidt2020catalog} 
\begin{equation}
4 \overset{\infty}{\underset{k =0}{\sum}}\; \frac{\left(-1\right)^{k}\,x^{2\,k +1}}{1-x^{2\,k +1}}
 = \vartheta_{3}\! \left(0, x\right)^{2}-1 \hspace{20pt} 0<x<1
\label{Ls}
\end{equation}
where $\vartheta_{3}\! \left(0, x\right)$ is the Jacobi theta function. Fundamentally, \eqref{Eq3} reduces to the elementary identity
\begin{equation}
\overset{\infty}{\underset{k =1}{\sum}}\; \frac{1}{1+a^{k}} = 
\,-\overset{\infty}{\underset{k =1}{\sum}}\; \frac{1}{1-a^{k}}+2 \overset{\infty}{\underset{k =1}{\sum}}\; \frac{1}{1-a^{2\,k}}\,.
\label{Triv}
\end{equation}


\subsection{Example \ref{sec:T1p4}} \label{sec:T1p4}
   We start by noting the identity 
   \begin{equation}
   (a^m+e^{-t})^{-1}-(a^m+e^t)^{-1}=\frac{2\sinh t}{ (a^{2m}+2a^m\cosh t+1)},
   \label{GId}
   \end{equation} 
 and, with \eqref{eqno(4)} in mind, utilize
\begin{equation}
f(a)=1/(a+x),~~ \quad F(s) = \pi\,a^{s - 1}/\sin(\pi\,s)
\label{FsEq4}
\end{equation}
to obtain
\begin{align} \nonumber
\overset{\infty}{\underset{n =0}{\sum}}\; \frac{1}{a^{n}+{\mathrm e}^{t}}
 & =
\frac{\left(2\,t +\ln \! \left(a \right)\right) {\mathrm e}^{-t}}{2\,\ln \! \left(a \right)}-\overset{\infty}{\underset{k =1}{\sum}}\; \frac{\left(-1\right)^{k}\,{\mathrm e}^{-\left(k +1\right) t}}{a^{k}-1}\\&+
\frac{2\,\pi \,}{\ln \! \left(a \right)}\Im \; \overset{\infty}{\underset{k =1}{\sum}}\; {\exp}\left({\frac{\left(2\,i\,\pi \,k -\ln \left(a \right)\right) t}{\ln \left(a \right)}}\right)\,\mathrm{csch}\! \left(\frac{2\,k\,\pi^{2}}{\ln \left(a \right)}\right)\,
\label{Pr}
\end{align} 
by evaluating the residues as before. Let $t:= -t$, subtract, and with $a:= a^m$, after comparing with \eqref{GId} we have
\begin{align} \nonumber
\overset{\infty}{\underset{n =0}{\sum}}\; \frac{1}{a^{2\,m\,n}+2\,a^{m\,n}\,\cosh \! \left(t \right)+1}
 = 
\frac{1}{2}\,-\frac{2\,\pi \,\coth \! \left(t \right) }{m\,\ln \! \left(a \right)}\overset{\infty}{\underset{k =1}{\sum}}\; \frac{\sin \left(\frac{2\,\pi \,k\,t}{m\,\ln \left(a \right)}\right)}{\sinh \left(\frac{2\,k\,\pi^{2}}{m\,\ln \left(a \right)}\right)}\\
-\frac{1}{\sinh \! \left(t \right)}\overset{\infty}{\underset{k =1}{\sum}}\! \frac{\left(-1\right)^{k}\,\sinh \left(\left(k +1\right) t \right)}{a^{m\,k}-1}-\frac{t\,\coth \! \left(t \right) }{m\,\ln \! \left(a \right)}\,.
\label{Ex4Ans}
\end{align}
By evaluating \eqref{Ex4Ans} in the limit $t\rightarrow 0$ and setting $a^m:=b$, we find

\begin{equation}
\overset{\infty}{\underset{n =0}{\sum}}\; \frac{1}{\left(b^{n}+1\right)^{2}}
 = 
\frac{1}{2}-\frac{1}{\ln \! \left(b \right)}-\frac{4\,\pi^{2}}{\ln^2 \! \left(b \right)}\overset{\infty}{\underset{k =1}{\sum}}\; k\,\mathrm{csch}\! \left(\frac{2\,k\,\pi^{2}}{\ln \left(b \right)}\right)+\overset{\infty}{\underset{k =2}{\sum}}\; \frac{k\,\left(-1\right)^{k}}{b^{k -1}-1}\,.
\label{Ct1}
\end{equation}
However, by expanding the denominator and transposing the resulting series (e.g. \eqref{Aik}), it is easy to write

\begin{equation}
\overset{\infty}{\underset{k =2}{\sum}}\; \frac{k\,\left(-1\right)^{k}}{b^{k -1}-1}
 = 
\overset{\infty}{\underset{n =1}{\sum}}\; \frac{1+2\,b^{n}}{\left(1+b^{n}\right)^{2}},\quad b>1,
\label{Sumbk}
\end{equation} 
so that \eqref{Ct1} reduces to a transformation between similar generalized Lambert series (\cite{AgRama} )

\begin{equation}
\overset{\infty}{\underset{n =1}{\sum}}\; \frac{b^{n}}{\left(1+b^{n}\right)^{2}}-\frac{1}{2\,\ln \! \left(b \right)}
 = 
-\frac{1}{8}+\frac{4\,\pi^{2}}{\ln \! \left(b \right)^{2}}\; \overset{\infty}{\underset{j =0}{\sum}}\; \frac{q^{1+2j}}{\left(q^{1+2j} -1\right)^{2}} \quad b>1,
\label{Ct1C}
\end{equation}
where
\begin{equation}
q = {\mathrm e}^{{2\,\pi^{2}}/{\ln \left(b \right)}}\,.
\label{Qid}
\end{equation}
{\bf Remarks:}
\begin{itemize}
\item{}
Because the sum on the right-hand side of \eqref{Ct1C} effectively vanishes exponentially as $b\rightarrow 1$, we obtain 
\begin{equation}
\underset{b \rightarrow 1}{\mathrm{lim}}\! \left(\overset{\infty}{\underset{j =1}{\sum}}\; \frac{b^{j}}{\left(1+b^{j}\right)^{2}}-\frac{1}{2\,\ln \! \left(b \right)}\right)
 = -{\frac{1}{8}}\,;
\label{Blim}
\end{equation}
\item{}
By setting $b=\exp(2\pi\,a),~ a>0$, we arrive at
\begin{equation}
\overset{\infty}{\underset{j =1}{\sum}}\; {\rm{sech}^{2}} \! \left(j\,\pi \,a \right)-\frac{1}{a^{2}}\;\overset{\infty}{\underset{j =0}{\sum}}\; {\rm{csch}}^{2} \left(\left(j +{1}/{2}\right)\frac{ \pi}{a}\right)
 = \frac{1}{\pi \,a}-\frac{1}{2}\,,
\label{CosSq}
\end{equation}
a known result when $a=1$, if one notes that the misprint $\Gamma(1/4)^2$, listed by both Hansen \cite[Eq. (43.8.12)]{Hansen} and Erdelyi et.al. \cite[Eq. 5.3.7(13)]{Erdelyi1}, should be $\Gamma(1/4)^4$. 

\item{}
In the case that $a=\exp(t/m)$, with $t>0$, \eqref{Ex4Ans} becomes

\begin{equation}
\overset{\infty}{\underset{n =0}{\sum}}\; \frac{1}{{\mathrm e}^{2\,n\,t}+2\,{\mathrm e}^{n\,t}\,\cosh \! \left(t \right)+1}
 = 
\frac{1}{2}-\frac{\cosh \! \left(t \right)}{\sinh \! \left(t \right)}-\frac{1}{\sinh \! \left(t \right)}\overset{\infty}{\underset{k =1}{\sum}}\; \frac{\left(-1\right)^{k}\,\sinh \left(\left(k +1\right) t \right)}{{\mathrm e}^{t\,k}-1}\,.
\label{Ex4ANs2}
\end{equation}
\end{itemize}

\subsection{Example \ref{sec:T1p5}} \label{sec:T1p5}

   
Here we consider the Mellin transform pair $f(x)=e^{-bx}$ and $F(s)=b^{-s}\Gamma(s)$ giving

\begin{equation}
\overset{\infty}{\underset{j =0}{\sum}}\; {\mathrm e}^{-b^{j}} = 
\frac{1}{2\,\pi\,i}\int_{c -i\,\infty }^{c +i\,\infty }\frac{\Gamma \left(s \right)}{1-b^{-s}}d s 
\label{Ein1}
\end{equation}
with $\Re(b)>1$ and $c>0$. Shifting the contour such that $-1<c<0$ produces
\begin{equation}
\overset{\infty}{\underset{j =0}{\sum}}\; {\mathrm e}^{-b^{j}} = \frac{1}{2}
-\frac{\gamma}{\ln \! \left(b \right)}+\frac{1}{2\,\pi\,i} \int_{c -i\,\infty }^{c +i\,\infty }\frac{\Gamma \left(s \right)}{1-b^{-s}}d s +\frac{2}{\ln \! \left(b \right)} \overset{\infty}{\underset{j =1}{\sum}}\; \Re \! \left(\Gamma \! \left(\frac{2\,i\,\pi \,j}{\ln \left(b \right)}\right)\right)
\label{E2}
\end{equation}
by taking into account the residues of the poles at $s=0$ and $s=2\pi\,i\,j/\ln(b)$. Further shifts of the contour $N$ units to the left following an obvious change of variables, yields

\begin{align} \nonumber
\overset{\infty}{\underset{j =0}{\sum}}\; {\mathrm e}^{-b^{j}} = &
\frac{1}{2}-\frac{\gamma}{\ln \! \left(b \right)}+\frac{2}{\ln \! \left(b \right)}\overset{\infty}{\underset{j =1}{\sum}}\;\Re \left(\Gamma \! \left(\frac{2\,i\,\pi \,j}{\ln \left(b \right)}\right)\right)\\+ &\left[ \frac{1}{2\,\pi}\int_{-\infty}^{\infty}\frac{\Gamma \left(c_{N}+i\,v \right)}{1-b^{-i\,v -c_{N}}}d v +\overset{N}{\underset{k =1}{\sum}}\; \frac{\left(-1\right)^{k}}{\Gamma \! \left(1+k \right) \left(1-b^{k}\right)}\right] 
\label{E3A}
\end{align} 
where again $b>1$ and $-N-1<c_{N}<-N$. Since the terms enclosed in brackets ([..]) contain the only $N$ dependence, if we consider the case that $N\rightarrow \infty$, the sum of the enclosed terms must remain constant, and since the sum clearly converges, so must the integral. Since the integral does not vary over the range $-N-1<c_{N}<-N$, this allows us to choose $c_{N}=-N-1/2$ and note first that
\begin{equation}
\Gamma \! \left(-N -\frac{1}{2}+i\,v \right) = 
\frac{\Gamma \! \left(-\frac{1}{2}+i\,v \right)}{\overset{N -1}{\underset{j =0}{\prod}}\! \left(-N +j -\frac{1}{2}+i\,v \right)}
\label{Gid}
\end{equation}
and second that
\begin{equation}
\underset{N \rightarrow \infty}{\mathrm{lim}}\; \frac{1}{\overset{N -1}{\underset{j =0}{\prod}}\;\left(-N +j -\frac{1}{2}\right)}
\sim {\left(\frac{e}{N}\right)^{N}}/{N}
\label{Prod}
\end{equation}
and therefore if the contour is moved such that  $c_{N}\rightarrow \,-\,\infty$, the integral vanishes, leaving

\begin{equation}
\overset{\infty}{\underset{j =0}{\sum}}\; {\mathrm e}^{-b^{j}} = 
\frac{1}{2}-\frac{\gamma}{\ln \! \left(b \right)}-\overset{\infty}{\underset{k =1}{\sum}}\; \frac{\left(-1\right)^{k}}{\Gamma \! \left(1+k \right) \left(b^{k}-1\right)}+\frac{2}{\ln \! \left(b \right)}\overset{\infty}{\underset{j =1}{\sum}}\; \Re\left(\Gamma \! \left(\frac{2\,i\,\pi \,j}{\ln \left(b \right)}\right)\right)\,.
\label{E4a}
\end{equation}
Further, since $b>1$, we can expand the denominator term in the first sum on the right-hand side of \eqref{E4a}, interchange the two sums and eventually identify
\begin{equation}
\overset{\infty}{\underset{k =1}{\sum}}\; \frac{\left(-1\right)^{k}}{\Gamma \! \left(1+k \right) \left(b^{k}-1\right)}
 = 
\overset{\infty}{\underset{j =1}{\sum}}\left({\mathrm e}^{-{1}/{b^{j}}}-1\right)\,,
\label{Sx}
\end{equation}
in which case we find 
\begin{equation}
\overset{\infty}{\underset{j =1}{\sum}}\; \Re \! \left(\Gamma \! \left(\frac{2\,i\,\pi \,j}{\ln \! \left(b \right)}\right)\right)
 = 
\frac{\ln \! \left(b \right)}{2} \,\overset{\infty}{\underset{j =0}{\sum}}\! \left({\mathrm e}^{-{1}/{b^{j}}}-1+{\mathrm e}^{-b^{j}}\right)+\left(\frac{1}{4}-\frac{1}{2\,e}\right) \ln \! \left(b \right)+\frac{\gamma}{2}\,,
\label{E4B}
\end{equation}
an identity that could also be rewritten as
\begin{equation}
\overset{\infty}{\underset{j =1}{\sum}}\; \Re \! \left(\Gamma \! \left(\frac{2\,i\,\pi \,j}{b}\right)\right)
 = 
\frac{b}{2} \overset{\infty}{\underset{j =0}{\sum}}\; \left({\mathrm e}^{-{1}/{{\mathrm e}^{j\,b}}}-1+{\mathrm e}^{-{\mathrm e}^{j\,b}}\right)+\left(\frac{1}{4}-\frac{1}{2e}\right) b +\frac{\gamma}{2}
\label{E4B1}
\end{equation}
by setting $b:=\exp(b)$.

\section{Examples based on Theorem 2} \label{sec:Dirichlet}

\subsection{Example \ref{sec:Ex2p1}} \label{sec:Ex2p1}


Continuing from the previous Section, consider the transform pair \eqref{eqno(7)}

\begin{align}
F(v)&=\Gamma(v)\zeta(v)\\
f(a)&=1/(\exp(a)-1)
\label{MellEx1}
\end{align}
leading to the identity

\begin{equation}
\frac{1}{2\,\pi}\int_{-\infty}^{\infty}\zeta \! \left(-s \left(c+i\,v  \right)\right) \Gamma \! \left(c+i\,v \right) \zeta \! \left(c+i\,v \right)d v = \overset{\infty}{\underset{j =1}{\sum}}\; \frac{1}{{\exp}({j^{{-s}})}-1}
\label{F11}
\end{equation}
after applying \eqref{F0}, where both sides converge if $s\in\Re\,,~s<0,~c>1~$ and $c>-1/s$. By shifting the contour left, variations arise by evaluating the appropriate residues as follows:

\begin{align} \nonumber
&\frac{1}{2\,\pi}\int_{-\infty}^{\infty}\zeta \! \left(-s \left(c+i\,v  \right)\right) \Gamma \! \left(c+i\,v  \right) \zeta \! \left(c+i\,v  \right)d v-\overset{\infty}{\underset{j =1}{\sum}}\; \frac{1}{\exp{(j^{^{-s}}})-1} \\&\hspace{45pt}= \left\{\begin{array} {lr}
0 &~~~~~ c>1\mbox{ and } c>-1/s,\\\\
\Gamma \! \left(-{1}/{s}\right) \zeta \! \left(-{1}/{s}\right)/s &~~~~~ c>1\mbox{ and } c<-1/s,\\\\
-\zeta \! \left(-s \right) &~~~~~ 0<c<1\mbox{ and } c>-1/s\\\\
\Gamma \! \left(-{1}/{s}\,\right) \zeta \! \left(-{1}/{s}\right)/s-\zeta \! \left(-s \right) &~~~~~ 0<c<1\mbox{ and } c<-1/s\\\\
{\Gamma \! \left(-{1}/{s}\right) \zeta \! \left(-{1}/{s}\right)}/{s}-\zeta \! \left(-s \right) 
-\overset{N}{\underset{j =0}{\sum}}\;{\displaystyle \frac{\left(-1\right)^{j}\,\zeta \! \left(s\,j \right) \zeta \! \left(-j \right)}{\Gamma \! \left(1+j \right)}} &~~~~~-N-1<c<-N\,,\\\\
\end{array}\right.
\label{Fall}
\end{align}
where $N=0,1,...$. 

We also now consider the transform pair \eqref{eqno(7a)} leading to

\begin{equation}
\frac{1}{2\,\pi}\int_{-\infty}^{\infty}\zeta \! \left(-s \left(c+i\,v  \right)\right) \Gamma \! \left(c+i\,v  \right) \zeta \! \left(c+i\,v  \right) \left(1-2^{1-i\,v -c }\right)d v
 = 
\overset{\infty}{\underset{j =1}{\sum}}\;\frac{1}{\exp({j^{-s})}+1}
\label{F1}
\end{equation}

where we require $s\in\Re$, $s<0$ and $c>-1/s$. Again, if the contour is shifted left, we find that various residues must be incorporated depending on the relative values of $s$ and $c$. Specifically

\begin{align} \nonumber
&\frac{1}{2\,\pi}\int_{-\infty}^{\infty}\zeta \! \left(-s \left(c+i\,v\ \right)\right) \Gamma \! \left(c+i\,v \right) \zeta \! \left(c+i\,v \right) \left(2^{1-v\,i-c }-1\right)d v+\overset{\infty}{\underset{j =1}{\sum}}\; \frac{1}{{\mathrm e}^{j^{^{- s}}}+1}
\\&\hspace{15pt}= 
 \left\{\begin{array} {lr}
0 &~~~~~~ c>-1/s,\\\\
\Gamma \! \left(-{1}/{s}\right) \zeta \! \left(-{1}/{s}\right) \left(-1+2^{1+{1}/{s}}\right)/s &~~~~~ c>0\mbox{ and } c<-1/s,
\\\\
\Gamma \! \left(-{1}/{s}\right) \zeta \! \left(-{1}/{s}\right) \left(-1+2^{1+{1}/{s}}\right)/s   \\\\
 -\overset{N}{\underset{j =1}{\sum}}\; \frac{\left(2^{2\,j}-1\right)}{\Gamma \! \left(2\,j +1\right)}\,\zeta \! \left(\left(2\,j -1\right) s \right) B_{2\,j} -{1}/{4} &~~~~~-N-1<c<-N\,,\\
\end{array}\right.
\label{Fall4}
\end{align}
where $B_{2j}$ are Bernoulli numbers (see \eqref{ZrefB}). In the case of equality, half the residue at that point must be included.

\subsubsection{The case $s=-1$} \label{sec:seqm1}
By taking the appropriate limit in \eqref{Fall}, let $s\rightarrow -1$, which, with $c=1/2$, gives

\begin{equation} 
\frac{1}{2\,\pi}\int_{-\infty}^{\infty}\zeta \! \left(\frac{1}{2}+i\,v \right)^{2}\,\Gamma \! \left(\frac{1}{2}+i\,v \right)d v
=\overset{\infty}{\underset{j =1}{\sum}}\; \frac{1}{{\mathrm e}^{j}-1}-\gamma\,,
\label{F5}
\end{equation}

which can be rewritten as

\begin{equation}
\frac{1}{\sqrt{\pi}}\int_{0}^{\infty}\frac{{| \zeta \left(\frac{1}{2}+i\,v \right)|}^{2}\,\cos \left(2\,\alpha \left(\frac{1}{2}+i\,v \right)+\theta \left(\frac{1}{2}+i\,v \right)\right)}{\sqrt{\cosh \left(\pi \,v \right)}}d v
=\overset{\infty}{\underset{j =1}{\sum}}\; \frac{1}{{\mathrm e}^{j}-1}-\gamma
\label{F5b}
\end{equation}

where
\begin{equation}
\zeta \left(\frac{1}{2}+i\,v \right)\equiv e^{i\,\alpha\left(1/2+iv\right)}| \zeta \left(\frac{1}{2}+i\,v \right)|
\label{Zdef}
\end{equation}

and 
\begin{equation}
\Gamma \left(\frac{1}{2}+i\,v \right)\equiv e^{i\,\theta\left(1/2+iv\right)}| \Gamma \left(\frac{1}{2}+i\,v \right)|
\label{Gdef}
\end{equation}

by writing both in polar form and employing the identity  \cite[Eq. 5.4.4]{NIST} 
\begin{equation}
\left| \Gamma \left(\frac{1}{2}+i\,v \right)\right|=\sqrt{\pi/\cosh{\pi\,v}}.
\label{GhalfId}
\end{equation}  

From the identity (\cite[Eq. (6.15)]{Milgram_Exploring})

\begin{equation}
\alpha \! \left({1}/{2}+i\,v \right) = 
\frac{v\,\ln \! \left(2\,\pi \right)}{2}-\frac{\theta \! \left(\frac{1}{2}+i\,v \right)}{2}-\frac{9\,\pi}{8}+\frac{\arctan \! \left({\mathrm e}^{\pi \,v}\right)}{2}\,,
\label{AId3}
\end{equation}
\eqref{F5b} then identifies

\begin{align} \nonumber
\int_{0}^{\infty}&\frac{{| \zeta \! \left(\frac{1}{2}+i\,v \right)|}^{2} }{\cosh \! \left(\pi \,v \right)} \left(\cos \! \left(v\,\ln \! \left(2\,\pi \right)\right) \cosh \! \left(\frac{\pi \,v}{2}\right)-\sin \! \left(v\,\ln \! \left(2\,\pi \right)\right) \sinh \! \left(\frac{\pi \,v}{2}\right)
\right) d v \\
 &= 
\sqrt{\pi}\left(\overset{\infty}{\underset{j =1}{\sum}}\; \frac{1}{{\mathrm e}^{j}-1}-\gamma\right)
\label{F5b4}
\end{align}
after applying elementary trigonometric identities and simplification. We now consider \eqref{Fall4} with the same limit $s\rightarrow -1$, and, comparing with \eqref{F5}, arrive at

\begin{equation}
\int_{-\infty}^{\infty}\zeta \! \left(\frac{1}{2}+i\,v \right)^{2}\,\Gamma \! \left(\frac{1}{2}+i\,v \right) 2^{-i\,v}d v
 = 
\sqrt{2}\,\pi  \left(2\; \overset{\infty}{\underset{j =1}{\sum}}\; \frac{1}{{\mathrm e}^{2\,j}-1}-\gamma +\ln \! \left(2\right)\right)\,.
\label{F4c}
\end{equation}
Applying the same identities as above, yields the equivalent form

\begin{align} \nonumber
\int_{-\infty}^{\infty}&\frac{{| \zeta \! \left(\frac{1}{2}+i\,v \right)|}^{2}}{\cosh \! \left(\pi \,v \right)}\left(\sinh \! \left(\frac{\pi \,v}{2}\right) \sin \! \left(v\,\ln \! \left(\pi \right)\right)-\cosh \! \left(\frac{\pi \,v}{2}\right) \cos \! \left(v\,\ln \! \left(\pi \right)\right)\right) d v \\ &
 = 
-\sqrt{2\,\pi} \left(2\; \overset{\infty}{\underset{j =1}{\sum}}\; \frac{1}{{\mathrm e}^{2\,j}-1}-\gamma +\ln \! \left(2\right)\right)\,,
\label{F4dR2}
\end{align}
a companion to \eqref{F5b4}.

\subsubsection{$s=-2n$} \label{sec:seq2n}

\begin{itemize}
\item{Case $c=0^{-}$}


Consider the case $s=-2n$ where $n=1,2,\dots$. In that eventuality, all poles in \eqref{Fall} corresponding to $c<0$ vanish as does the finite sum with $N\geq1$, so the original contour can be moved with impunity as far to the left (where it does NOT vanish) as one wishes. 
Of more interest, with $c=0^{-}$, since the pole at $v=0$ is imaginary, by adding $\pi/4=$ half the residue at $v=0$, \eqref{Fall} becomes

\begin{align} \nonumber
\int_{-\infty}^{\infty}\Re \! \left(\zeta \! \left(2\,i\,v\,n \right) \right.&\left. \Gamma \! \left(i\,v \right) \zeta \! \left(i\,v \right)\right)d v
 \\ &
 =2\,\pi  \overset{\infty}{\underset{j =1}{\sum}}\; \frac{1}{{\mathrm e}^{j^{2\,n}}-1}-\frac{\pi}{n} \,\zeta \! \left(\frac{1}{2\,n}\right) \Gamma \! \left(\frac{1}{2\,n}\right)-2\,\zeta \! \left(2\,n \right) \pi -\frac{\pi}{4}
\label{Fs2nc0}
\end{align}
and the sum of of \eqref{Fall} and \eqref{Fall4} becomes
\begin{align} \nonumber
\int_{-\infty}^{\infty}\Re \! \left(2^{-i\,v}\right.&\left.\,\zeta \! \left(2\,i\,v \,n\right) \Gamma \! \left(i\,v \right) \zeta \! \left(i\,v \right)\right)d v
\\ & = 
2\,\pi  \overset{\infty}{\underset{j =1}{\sum}}\; \frac{1}{{\mathrm e}^{2j^{2\,n}}-1}-\frac{\pi \, 2^{-{1}/{2\,n}}}{n}\,\zeta \! \left(\frac{1}{2\,n}\right)\Gamma \! \left(\frac{1}{2\,n}\right)-\zeta \! \left(2\,n \right) \pi -\frac{\pi}{4}\,.
\label{Fb1c0}
\end{align}

If we now define the Riemann function
\begin{equation}
\Upsilon(s,b)\equiv \zeta(s)\,b^{-s/2}\,\Gamma(s/2)
\label{UpsDef}
\end{equation}
which is well-known to satisfy

\begin{equation}
\Upsilon(s,\pi)=\Upsilon(1-s,\pi)
\label{UpsRev}
\end{equation}
due to the reflection property of $\Gamma(s/2)$ and the functional equation of $\zeta(s)$, with $n=1$, \eqref{Fb1c0} can be rewritten

\begin{equation}
\int_{-\infty}^{\infty}\Re \! \left(\Upsilon \! \left(2\,i\,v , 2\right) \zeta \! \left(i\,v \right)\right)d v
 = 
-\frac{\zeta \! \left(\frac{1}{2}\right) \pi^{\frac{3}{2}}\,\sqrt{2}}{2}-\frac{\pi^{3}}{6}+2\,\pi  \overset{\infty}{\underset{j =1}{\sum}}\; \frac{1}{{\mathrm e}^{2j^{{2}}}-1}-\frac{\pi}{4}.
\label{Fb3Pr}
\end{equation}

From \eqref{UpsDef} and\eqref{UpsRev} we have the more general form of \eqref{UpsRev}

\begin{equation}
\Upsilon(s,b) = \Upsilon(1 - s,b)\left(\frac{b}{\pi}\right)^{(1/2 - s)}
\label{UpsRev2}
\end{equation}
so that taking \eqref{UpsRev2} into account yields a very different form of \eqref{Fb1c0}, that being

\begin{align} \nonumber
\int_{-\infty}^{\infty}\Re \!& \left(\zeta \! \left(1-2\,i\,v \right) \left(\frac{\pi^{2}}{2}\right)^{i\,v}\,\Gamma \! \left(\frac{1}{2}-i\,v \right) \zeta \! \left(i\,v \right)\right)d v \\
& = 
\pi^{\frac{3}{2}} \left(-\zeta \! \left(\frac{1}{2}\right)\sqrt{\frac{\pi}{2}}-\frac{\pi^{2}}{6}+2 \overset{\infty}{\underset{j =1}{\sum}}\; \frac{1}{{\mathrm e}^{2\,j^{2}}-1}-\frac{1}{4}\right)\,.
\label{Fb3aPr}
\end{align}

{\bf Remark:} Since $|\Gamma \! \left(i\,v \right)|^2={\pi}/{(v\, \sinh(\pi\,v))}$, all integrals are convergent.

In the case that $n\rightarrow \infty$, \eqref{Fs2nc0} and \eqref{Fb1c0} respectively become

\begin{equation}
\underset{n \rightarrow \infty}{\mathrm{lim}}\! \int_{-\infty}^{\infty}\Re \! \left(\zeta \! \left(2\,i\,v\,n \right) \Gamma \! \left(i\,v \right) \zeta \! \left(i\,v \right)\right)d v 
 = \frac{\pi  \left(13-5\,{\mathrm e}\right)}{4\,{\mathrm e}-4}
\label{Fs2nc0Lim}
\end{equation}
and
\begin{equation}
\underset{n \rightarrow \infty}{\mathrm{lim}}\int_{-\infty}^{\infty}\Re \! \left(2^{-i\,v}\,\zeta \! \left(2\,i\,v\,n \right) \Gamma \! \left(i\,v \right) \zeta \! \left(i\,v\right)\right)d v 
 = 
\frac{\pi  \left(9-{\mathrm e}^{2}\right)}{4 \left({\mathrm e}^{2}-1\right)}\,.
\label{Fb1c0Lim}
\end{equation}

\item{Case $c=1/(2n)$}

Other interesting cases arise when $c=1/(2n)=-1/s$ when the poles of the integrand appear only in the imaginary part and we must use the corresponding half residues in \eqref{Fall} and \eqref{Fall4}. Specifically, and of interest, the following are obtained by combining the two cases appropriately, to yield

\begin{align} \nonumber
\int_{-\infty}^{\infty}&\Re \! \left(\zeta \! \left(1+2\,i\,v\,n \right) \Gamma \! \left(\frac{1}{2\,n}+i\,v \right) \zeta \! \left(\frac{1}{2\,n}+i\,v \right)\right)d v \\
& = 
2\,\pi  \overset{\infty}{\underset{j =1}{\sum}}\; \frac{1}{{\mathrm e}^{j^{2\,n}}-1}-\pi \,\zeta \! \left(\frac{1}{2\,n}\right) \Gamma \! \left(1+\frac{1}{2\,n}\right)-2\,\pi\,\zeta \! \left(2\,n \right) 
\label{Ga1}
\end{align}
and
\begin{align} \nonumber
&\int_{-\infty}^{\infty}\Re \! \left(\left(2^{-{1}/{(2\,n)}-i\,v}-1\right) \zeta \! \left(1+2\,i\,v\,n\right) \zeta \! \left(\frac{1}{2\,n}+i\,v \right) \Gamma \! \left(\frac{1}{2\,n}+i\,v \right)\right)d v=\pi\,\zeta \! \left(2\,n \right)\\
 &-2\,\pi\left( 
 \overset{\infty}{\underset{j =1}{\sum}}\; \frac{1}{{\mathrm e}^{j^{2\,n}}-1}- \overset{\infty}{\underset{j =1}{\sum}}\; \frac{1}{{\mathrm e}^{2j^{2\,n}}-1}\right)+\pi \left(1-2^{-\frac{1}{2\,n}}\right) \zeta \! \left(\frac{1}{2\,n}\right) \Gamma \! \left(1+\frac{1}{2\,n}\right)\,. 
\label{Gba1}
\end{align}
In contrast to \eqref{Ga1}, in the case that $n\rightarrow\infty$, the integrand \eqref{Gba1} converges at $v=0$, and we find
\begin{equation}
\underset{n \rightarrow \infty}{\mathrm{lim}}\! \int_{-\infty}^{\infty}\Re \! \left(\left(2^{-i\,v}-1\right) \zeta \! \left(1+2\,i\,v\,n \right) \zeta \! \left(i\,v \right) \Gamma \! \left(i\,v \right)\right)d v 
 = 
{\displaystyle\pi}\,\frac{\left({\mathrm e}^{2}-2\,{\mathrm e}-1\right) }{\left({\mathrm e}-1\right) \left({\mathrm e}+1\right)}.
\label{Gb3}
\end{equation}

\item{Case $c=1/n$ and $c=1/(4n)$}

Other special cases abound, among which we consider $c=1/n$ and $c=1/(4n)$ to respectively yield
\begin{align} \nonumber
\int_{-\infty}^{\infty}&\Re \! \left(\zeta \! \left( 2+2\,i\,v\,n\right) \Gamma \! \left(\frac{1}{n}+i\,v \right) \zeta \! \left(\frac{1}{n}+i\,v \right)\right)d v \\ &
 = 
2\,\pi  \overset{\infty}{\underset{j =1}{\sum}}\; \frac{1}{{\mathrm e}^{j^{2\,n}}-1}-2\,\pi\,\zeta \! \left(2\,n \right) 
\label{GA1}
\end{align}
\begin{align} \nonumber
\int_{-\infty}^{\infty}&\Re \! \left(\zeta \! \left(\frac{1}{2}+2\,i\,v\,n \right) \Gamma \! \left(\frac{1}{4\,n}+i\,v \right) \zeta \! \left(\frac{1}{4\,n}+i\,v \right)\right)d v \\ &
 = 
2\,\pi \, \overset{\infty}{\underset{j =1}{\sum}}\; \frac{1}{{\mathrm e}^{j^{2\,n}}-1}-\frac{\pi}{n} \,\zeta \! \left(\frac{1}{2\,n}\right) \Gamma \! \left(\frac{1}{2\,n}\right)-2\pi\,\zeta \! \left(2\,n \right) .
\label{G1A}
\end{align}
\end{itemize}

\subsubsection{Other values of $s$} \label{sec:Other}

For other values of $s$, interesting cases also arise. For example, if $s=-n$ and $c=1/(2n)$ we arrive at

\begin{align} \nonumber
\int_{-\infty}^{\infty}&\Re \! \left(\zeta \! \left(\frac{1}{2}+i\,v\,n \right) \Gamma \! \left(\frac{1}{2\,n}+i\,v \right) \zeta \! \left(\frac{1}{2\,n}+i\,v \right)\right)d v \\ &
 = 
2\,\pi \overset{\infty}{\underset{j =1}{\sum}}\; \frac{1}{{\mathrm e}^{j^{n}}-1}-2 \pi\,\zeta \! \left(\frac{1}{n}\right) \,\Gamma \! \left(1+\frac{1}{n}\right)-2\,\pi \,\zeta \! \left(n \right),
\label{Gd}
\end{align}
reducing to \eqref{F5} if $n=1$ and \eqref{G1A} if $n:=2n$. If $s=-1/n$ with $c=1/2$ we find

\begin{align} \nonumber
\int_{-\infty}^{\infty}&\Re \! \left(\zeta \! \left(\frac{1}{2\,n}+\frac{i\,v}{n}\right) \Gamma \! \left(\frac{1}{2}+i\,v \right) \zeta \! \left(\frac{1}{2}+i\,v \right)\right)d v \\&
 = 
2\,\pi  \overset{\infty}{\underset{j =1}{\sum}}\; \frac{1}{{\mathrm e}^{j^{{1}/{n}}}-1}-2\,\pi \,\zeta \! \left(n \right) \Gamma \! \left(n +1\right)-2\,\pi \,\zeta \! \left(\frac{1}{n}\right)
\label{Ge}
\end{align}
and, as $n\rightarrow\infty$ the first factor in the integrand reduces to $\zeta(0)=-1/2$ because most of the integrand originates near $v=0$ due to \eqref{GhalfId}, giving 
 
\begin{align} \nonumber
\int_{-\infty}^{\infty}&\Re \! \left(\Gamma \! \left(\frac{1}{2}+i\,v \right) \zeta \! \left(\frac{1}{2}+i\,v \right)\right)d v\\& = 
4\,\pi \underset{n \rightarrow \infty}{\mathrm{lim}}\! \left(\zeta \! \left(n \right) \Gamma \! \left(n +1\right)+\zeta \! \left(\frac{1}{n}\right)-\overset{\infty}{\underset{j =1}{\sum}}\; \frac{1}{{\mathrm e}^{j^{{1}/{n}}}-1}\right)
\,,
\label{GeAsy}
\end{align}
an identity whose right-hand side confounds numerical verification - however see \eqref{T2bIdX} below.

\subsection{Example \ref{sec:X2p2}} \label{sec:X2p2} 

Here we again consider the Mellin transform pair $f(x)=e^{-bx}$ and $F(s)=b^{-s}\Gamma(s)$ as in Section \ref{sec:T1p5}, giving, with $b>0$, 

\begin{equation}
\frac{1}{2\,\pi\,i} \int_{c -i\,\infty }^{c +i\,\infty }\zeta \! \left(-s\,v \right) b^{-v}\,\Gamma \! \left(v \right)d v=\overset{\infty}{\underset{j =1}{\sum}}\; {\mathrm e}^{-b\,j^{-s}} =\omega(b,-s)
\label{T2}
\end{equation}
valid for $-1<s<0$ and $c>-1/s$. After a change of variables, \eqref{T2} can also be written as
\begin{equation}
\frac{1}{2\,\pi}\int_{-\infty}^{\infty}\zeta \! \left(-s \left(c+i\,v  \right)\right) b^{-c-i\,v }\,\Gamma \! \left(c+i\,v  \right)d v=\overset{\infty}{\underset{j =1}{\sum}}\; {\mathrm e}^{-b\,j^{-s}}\,.
\label{T2x}
\end{equation}
Furthermore, since \eqref{T2} is an inverse Mellin transform, by inverting, {\bf iff $s<0$ and $b>0$}, we find


\begin{equation}
\zeta \! \left(-b\,s \right) \Gamma \! \left(b \right) = 
\int_{0}^{\infty}x^{b -1}\overset{\infty}{\underset{j =1}{\sum}}\; {\mathrm e}^{-x\,j^{-s}}d x,~~s<0,
\label{Jb}
\end{equation}
reducing to the classic results \cite[Eq. (2.4.1)]{Titch2} and \cite[Eq. (2.6.2)]{Titch2} if $s=-1$ (see \eqref{Sexp}) and $s=-2$ respectively.

\subsubsection{Case $0<c<-1/s$}

In the case that $0<c<-1/s$, which allows $s<-1$, we must include a residue term, so \eqref{T2x} becomes

\begin{equation}
\frac{1}{2\,\pi}\int_{-\infty}^{\infty}\zeta \! \left(-s \left(c+i\,v  \right)\right) b^{-c-i\,v }\,\Gamma \! \left(c+i\,v  \right)d v=\omega(b,-s) - b^{1/s}\,\Gamma(1-1/s)\,,
\label{T2xa}
\end{equation}
and further, if $c<0$, by moving the contour $N+1$ units to the left, \eqref{T2x} becomes

\begin{align} \nonumber
&\int_{-\infty}^{\infty}\zeta \! \left(-s \left(c_{N}+i\,v \right)\right) b^{-c_{N}-i\,v }\,\Gamma \! \left(c_{N}+i\,v \right)d v\\
 &= 2\,\pi\left( \omega(b,-s)
- \overset{N}{\underset{j =0}{\sum}}\; \frac{\zeta \! \left(s\,j \right) \left(-b \right)^{j}}{\Gamma \! \left(1+j \right)} -\,b^{{1}/{s}}\Gamma \! \left(1-{1}/{s}\right)\right),
\label{T2c}
\end{align}
where $c$ has been replaced by $c_{N}$ such that $-N-1<c_{N}<-N$, $N\in\mathcal{Z}$ and always $s<0.\\$

{\bf Case $s=-1\\$} \label{ssec:Seqm1}

In the case that $s=-1$, we find
\begin{equation}
\int_{-\infty}^{\infty}\zeta \! \left(c+i\,v  \right) b^{ -c-i\,v}\,\Gamma \! \left(c+i\,v  \right)d v
 = 2\,\pi  \left(\frac{1}{{\mathrm e}^{b}-1}-\frac{X}{b}\right)
\label{T2b1}
\end{equation}
where $X=0$ if $c>1$ and $X=1$ if $0<c<1$. In the case that $c=1$, the singularity of the integrand only occurs in the imaginary part, and, with $X=1/2$ we obtain
\begin{equation}
\int_{-\infty}^{\infty}\Re \! \left(\zeta \! \left(1+i\,v \right) b^{-i\,v }\,\Gamma \! \left(1+i\,v \right)\right)d v
 = 2\,\pi\,b \left(\frac{1}{{\mathrm e}^{b}-1}-\frac{1}{2\,b}\right)\,. 
\label{T2b2}
\end{equation}

\begin{figure}[h] 
\centering
\includegraphics[width=.7\textwidth]{{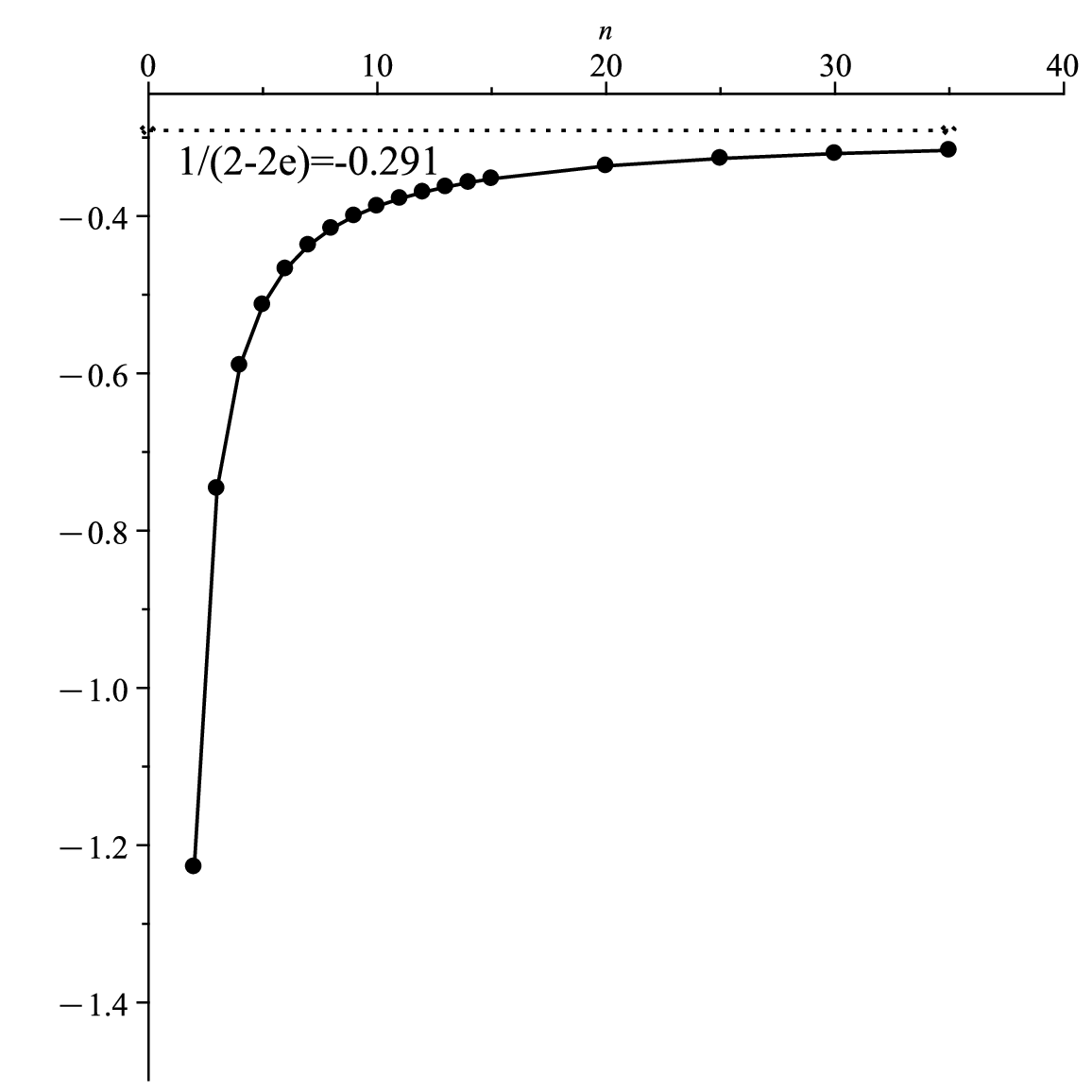}}
\caption{Approximations to the left-hand side of \eqref{T2bIdX} for increasing values of $n$ compared to the numerical value of the right-hand side. For $n>15$ an extraordinarily large number of digits was required.}
\label{fig:LimFig1}
\end{figure}

Similarly, in the case that $c=0$, we find, in exactly the same way

\begin{equation}
\int_{0}^{\infty}\Re \! \left(\zeta \! \left(i\,s\,v \right) b^{i\,v}\,\Gamma \! \left(-i\,v \right)\right)d v
 = 
\pi\, \omega(b,-s)-\pi\,b^{{1}/{s}}\,\Gamma \! \left(1-\frac{1}{s}\right)  +\frac{\pi}{4}\,,
\label{T2B}
\end{equation}
so that, if $s=-1$ we obtain

\begin{equation}
\int_{-\infty}^{\infty}\zeta \! \left(i\,v \right) b^{-i\,v}\,\Gamma \! \left(i\,v \right)d v
 = 2\,\pi  \left(\frac{1}{{\mathrm e}^{b}-1}-\frac{1}{b}+\frac{1}{4}\right)\,,
\label{V3sm1}
\end{equation}

and, if $c=1/2$, $b=1$ we find
\begin{equation}
\int_{-\infty}^{\infty}\zeta \! \left(\frac{1}{2}+i\,v \right) \Gamma \! \left(\frac{1}{2}+i\,v \right)d v
 = 2\,\pi \left(\frac{1}{{\mathrm e}-1}-1\right) \,.
\label{T2b3b}
\end{equation}
{\bf Remark:} Comparison of the right-hand sides of \eqref{GeAsy} and \eqref{T2b3b} identifies
\begin{equation}
\underset{n \rightarrow \infty}{\mathrm{lim}}\! \left(\overset{\infty}{\underset{j =1}{\sum}}\; \frac{1}{{\mathrm e}^{j^{{1}/{n}}}-1}-\zeta \! \left(n \right) \Gamma \! \left(n +1\right)\right)
 = \frac{1}{2-2\,{\mathrm e}}
\label{T2bIdX}
\end{equation}
or, alternatively, for large values of $n$,
\begin{equation}
    \overset{\infty}{\underset{j =1}{\sum}}\; \frac{1}{{\mathrm e}^{j^{{1}/{n}}}-1}
 \sim \sqrt{2\,\pi}\,n^{n +\frac{1}{2}}\,{\mathrm e}^{-n}\,.
\label{T2bAsy}
\end{equation}
This result can, with some difficulty, be tested numerically - see Figure \ref{fig:LimFig1}.\newline

{\bf Case $s=-2$ and yet another proof of the Poisson-Jacobi transform.}\newline

Consider the case $s=-2$ with $0<c<1/2$ yielding the identity

\begin{equation}
\int_{-\infty}^{\infty}\zeta \! \left(2\,c -2\,i\,v \right) b^{-c+i\,v }\,\Gamma \! \left(c-i\,v  \right)d v
 = 
2\,\pi  \overset{\infty}{\underset{j =1}{\sum}}\; {\mathrm e}^{-b\,j^{2}}-\frac{\pi^{\frac{3}{2}}}{\sqrt{b}}
\label{T3B}
\end{equation}
and, from \eqref{UpsDef} and \eqref{UpsRev2} the left-hand side of \eqref{T3B} satisfies

\begin{align} \nonumber
\int_{-\infty}^{\infty}\zeta \;& \left(2\,c -2\,i\,v \right) b^{-c+i\,v }\,\Gamma \! \left(c-i\,v  \right)d v\\
 &= 
b^{-c}\,\pi^{-\frac{1}{2}+2\,c}\int_{-\infty}^{\infty}\zeta \! \left(1-2\,c+2\,i\,v \right) \Gamma \! \left(\frac{1}{2}-c +i\,v \right) \left(\frac{\pi^{2}}{b}\right)^{-i\,v}d v\,.  
\label{T3Int}
\end{align}

Now, consider a reflection of integration variables, replacements $c:=1/2-c$ and $b:=\pi^2/b$, all the while retaining $0<c<1/2$, and find that \eqref{T3B} becomes

\begin{align} \nonumber
\int_{-\infty}^{\infty}\zeta \! \left(1-2\,c +2\,i\,v \right) \Gamma \! \left( \frac{1}{2}-c+i\,v \right) \left(\frac{\pi^{2}}{b}\right)^{-i\,v}d v \\
 = 
2\,\pi^{2-2\,c}\,b^{c -\frac{1}{2}}\overset{\infty}{\underset{j =1}{\sum}}\; {\mathrm e}^{-{\pi^{2}\,j^{2}}/{b}}-\pi^{{3}/{2}-2\,c}\,b^{c}\,.
\label{T3B1}
\end{align}

Comparison of \eqref{T3B}, \eqref{T3Int} and \eqref{T3B1} eventually leads to

\begin{equation}
\overset{\infty}{\underset{j =1}{\sum}}\; {\mathrm e}^{-\pi \,b\,j^{2}}-\sqrt{\frac{1}{b}}\; \overset{\infty}{\underset{j =1}{\sum}}\; {\mathrm e}^{-{\pi \,j^{2}}/{b}}
 = \frac{1}{2}\sqrt{\frac{1}{b}}-\frac{1}{2},
\label{SxId2}
\end{equation}

equivalent to the well-known Poisson-Jacobi transform \cite[page 124]{WW}
\begin{equation}
\theta_{3} \! \left(0,b \right) = 
\sqrt{\frac{1}{b}}\,\theta_{3} \! \left(0,\frac{1}{b}\right)\,.
\label{PJt}
\end{equation}

\subsubsection{Case $s=-2n$}


In  the case that $s=-2n$, the infinite sum in \eqref{T2c} vanishes except for the term corresponding to $j=0$, leading to

\begin{equation}
\int_{-\infty}^{\infty}\Re \! \left(\zeta \! \left(2\,n \left(c+i\,v  \right)\right) b^{-i\,v -c}\,\Gamma \! \left(i\,v +c \right)\right)d v
 = 
\pi\left( \,X +2\, \omega(b,2n)-2\,Y\,b^{-{1}/{(2n)}}\,\Gamma \! \left(1+\frac{1}{2\,n}\right)\right) 
\label{T2nA}
\end{equation}
where $\left\{X=1,~1/2,~0\right\}$ if $\left\{c<0,~c=0,~c>0\right\}$ and $\left\{Y=1,~1/2,~0\right\}~$ if $\\$$~\left\{c<1/(2n),~c=1/(2n),~c>1/(2n)\right\}$ respectively. A family of interesting integrals arises if we let $c=p/n$, where $p\in\Re$, yielding the following:
\begin{itemize}
\item{if $p<0$}
\begin{equation}
\int_{-\infty}^{\infty}\zeta \! \left(2\,p +2\,i\,v\,n \right) b^{-i\,v}\,\Gamma \! \left(\frac{p}{n}+i\,v \right)d v
 = 
\pi \, b^{{p}/{n}} \left(1+2\omega(b,2n)-2\,b^{-{1}/{(2\,n)}}\,\Gamma \! \left(1+\frac{1}{2\,n}\right)\right)\,;
\label{T3nA}
\end{equation}
\item{if $p>1/2$}
\begin{equation}
\int_{-\infty}^{\infty}\Re \! \left(\zeta \! \left( 2\,p+2\,i\,v\,n \right) b^{-i\,v}\,\Gamma \! \left(\frac{p}{n}+i\,v \right)\right)d v
 = 
2\,\pi \,b^{{p}/{n}} \omega(b,2n)\,;
\label{T3nB}
\end{equation}
\item{if $p=1/2$}
\begin{equation}
\int_{-\infty}^{\infty}\Re \! \left(\zeta \! \left(1+2\,i\,v\,n \right) b^{-i\,v}\,\Gamma \! \left(\frac{1}{2\,n}+i\,v \right)\right)d v
 = 
\pi  \left(2\,b^{{1}/{(2\,n)}} \omega(b,2n)-\Gamma \! \left(1+\frac{1}{2\,n}\right)\right)
\label{T3nC}
\end{equation}
\item{if $0<p<1/2$}
\begin{equation}
\int_{-\infty}^{\infty}\Re \! \left(\zeta \! \left(2\,p +2\,i\,v\,n \right) b^{-i\,v}\,
\Gamma \! \left(\frac{p}{n}+i\,v \right)\right)d v
 =  
2\,\pi \,b^{{p}/{n}}\left(\omega(b,2n)-b^{-{1}/{(2\,n)}}\,\Gamma \! \left(1+\frac{1}{2\,n}\right)\right)
\label{T3nD}
\end{equation}

\item{if $p=0$}
\begin{equation}
\int_{-\infty}^{\infty}\Re \! \left(\zeta \! \left(2\,i\,v\,n \right) b^{-i\,v}\,\Gamma \! \left(i\,v \right)\right)d v
 = 
\pi\left(\frac{1}{2}+2 \omega(b,2n)-2\,b^{-{1}/{(2\,n)}}\,\Gamma \! \left(1+\frac{1}{2\,n}\right)\right) \,.
\label{T3nE}
\end{equation}
\end{itemize}

{\bf Remarks:} \begin{itemize}
\item { The case \eqref{T3nD} covers the interior of the critical strip.}
\item { The above resolves a special case discussed in \cite[page 4]{KFGT}.}
\item {In any of the above, if $n=1$, we have the known \cite[Eq. (3)]{Romik} identity $\omega(\pi,2)=\frac{\pi^{1/4}}{2\Gamma(3/4)}-\frac{1}{2}$  .}
\item {The case \eqref{T3nA} reduces to the known result \cite[Eq. (3.9)] {Milgram2019} if we set $b=\pi$ and $n=1$.}
\item {Setting $\Gamma \! \left({p}/{n}+i\,v \right)\approx \Gamma \! \left(i\,v \right)$ will obtain valid numerical approximations for the case of large $n$, but any attempt to equate them at the limit $n\rightarrow \infty$ is incorrect, because that degenerates into the case $p=0$ and $p$ and $n$ are independent variables.}
\end{itemize}

\subsubsection{Case s=-1/n and $c<0$}


Here we consider the case $s=-1/n$, with $c<0$, allowing us to choose $-n-1<c<-n$, so from \eqref{T2c} and without loss of generality, let $c=-n-1/2$, leading to

\begin{align} \nonumber
\int_{-\infty}^{\infty}\Re \! &\left(\zeta \! \left(-1-\frac{1}{2\,n}+\frac{i\,v}{n}\right) b^{n +\frac{1}{2}-i\,v}\,\Gamma \! \left(-n -\frac{1}{2}+i\,v \right)\right)d v +2\,\pi  \overset{n}{\underset{j =0}{\sum}}\; \frac{\zeta \! \left(-\frac{j}{n}\right) \left(-b \right)^{j}}{\Gamma \! \left(1+j \right)} \\
 &= 
2\,\pi \overset{\infty}{\underset{j =1}{\sum}}\; {\mathrm e}^{-b\,j^{{1}/{n}}}-2\,\pi\,b^{-n}\,\Gamma \! \left(1+n \right)\,. 
\label{Tc1a}
\end{align}
If we now consider the limiting case $n\rightarrow \infty$, it is easy to discover that the integration term vanishes by writing

\begin{equation}
\Gamma \! \left(-n -\frac{1}{2}+i\,v \right) = 
-\frac{\pi  \left(-1\right)^{n}}{\overset{n}{\underset{j =0}{\prod}}\; \left(\frac{1}{2}+n -j -i\,v \right) \Gamma \! \left(\frac{1}{2}-i\,v \right) \cosh \! \left(\pi \,v \right)},
\label{GidPr}
\end{equation}

and the left-side summation term becomes

\begin{equation}
\underset{n \rightarrow \infty}{\mathrm{lim}}\; \overset{n}{\underset{j =0}{\sum}}\; \frac{\zeta \! \left(-\frac{j}{n}\right) \left(-b \right)^{j}}{\Gamma \! \left(1+j \right)}
 = -\frac{1}{2\,{\mathrm e}^{b}}
\label{SbId}
\end{equation}
because as $n\rightarrow \infty$, $\zeta \! \left(-\frac{j}{n}\right)\approx \zeta(0)$. The identity \eqref{SbId} can be numerically verified for a large range of the variable $b>0$. Therefore, we find

\begin{equation}
\underset{n \rightarrow \infty}{\mathrm{lim}}\; \left(\overset{\infty}{\underset{j =1}{\sum}}\; {\mathrm e}^{-b\,j^{{1}/{n}}}-b^{-n}\,\Gamma \! \left(1+n \right)\right)
 = -\frac{1}{2\,{\mathrm e}^{b}}\,,
\label{Sgenf}
\end{equation}
or equivalently, as $n\rightarrow \infty$,

\begin{equation}
\overset{\infty}{\underset{j =1}{\sum}}\; {\mathrm e}^{-b\,j^{\frac{1}{n}}}
\sim \sqrt{2\,\pi \,n} \left(\frac{n}{{\mathrm e}\,b}\right)^{n}\,.
\label{As2}
\end{equation}

The identity \eqref{Sgenf} can, as with \eqref{T2bIdX}, be tested numerically - see Figure \ref{fig:LimFig2} and the remark following \eqref{Q2} below. See also \cite[Eq. 2.5]{KFGT}.

 
\begin{figure}[h]
\centering
\includegraphics[width=.8\textwidth]{{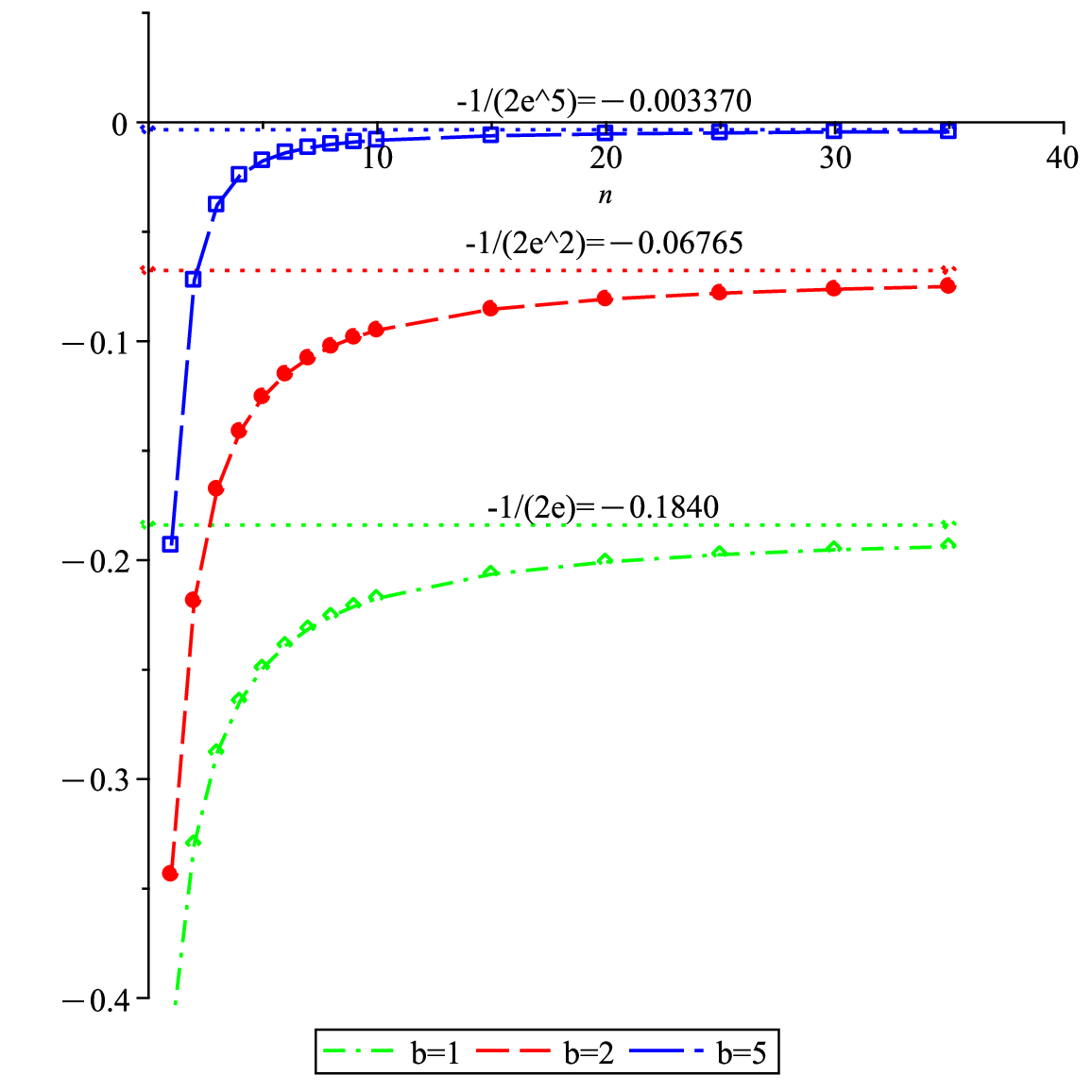}}
\caption{Approximations to the left-hand side of \eqref{Sgenf} for increasing values of $n$ compared to the numerical value of the right-hand side for three values of the parameter $b$.}
\label{fig:LimFig2}
\end{figure}

\subsubsection{Case $s=-2n$ and $c<0$}

In this case, again we choose $c=-n/2$ so \eqref{T2c} becomes

\begin{align} \nonumber
b^{n +\frac{1}{2}}\int_{-\infty}^{\infty}\zeta \! \left(2\,i\,n\,v -2\,n^{2}-n \right) b^{-i\,v}\,\Gamma \! \left(-n -\frac{1}{2}+i\,v \right)d v \\
 = 
\pi +2\,\pi  \left(\overset{\infty}{\underset{j =1}{\sum}}\; {\mathrm e}^{-b\,j^{2\,n}}-b^{-\frac{1}{2\,n}}\,\Gamma \! \left(1+\frac{1}{2\,n}\right)\right)
\label{Q1}
\end{align}
because only the term indexed by $j=0$ in the infinite sum included in \eqref{T2c} does not vanish. We now consider the limiting case $n=0$, leading to

\begin{equation}
\frac{1}{2\,\pi}\,\int_{-\infty}^{\infty}b^{-i\,v}\,\Gamma \! \left(-\frac{1}{2}+i\,v \right)d v
 = b^{-1/2}\,\left({1}/{{\mathrm e}^{b}}-1\right)
\label{Q2}
\end{equation}
by identifying $2n:=1/n$ in \eqref{Sgenf}.\newline 

{\bf Remark:} By a simple change of integration variables in \eqref{Q2}, wrapping the contour about the negative real axis and evaluating the residues so enclosed we arrive at an equivalent representation of \eqref{Q2}

\begin{equation}
\frac{1  }{2\,\pi\,i}\int_{-\,i\,\infty }^{\,i\,\infty }b^{-t}\,\Gamma \! \left(-\frac{1}{2}+t \right)d t
 = 
b^{-1/2}\,\overset{\infty}{\underset{n =1}{\sum}}\; \frac{b^{n} \left(-1\right)^{n}}{\Gamma \left(1+n \right)},
\label{Q2c}
\end{equation}
recognizing that the equivalence of the right-hand sides of \eqref{Q2} and \eqref{Q2c} is an elementary identity, yielding a  validity check of \eqref{Sgenf} since it is the basis for \eqref{Q2}.\newline

\subsubsection{Case $s=-(2n+1)$ and $c<0$}
As before, without loss of generality, we choose $c=-N/2,~N>0$ in \eqref{T2c}, whose general form becomes

\begin{align} \nonumber
&\int_{-\infty}^{\infty}\zeta \! \left({\left(n +1/2\right) \left(2\,i\,v -2\,N -1\right)}\right) b^{N +\frac{1}{2}-i\,v}\,\Gamma \! \left(-N -{1}/{2}+i\,v \right)d v \\
&
+2\,\pi  \overset{N}{\underset{j =0}{\sum}}\; \frac{\zeta \; \left(\left(-2\,n -1\right) j \right) \left(-b \right)^{j}}{\Gamma \! \left(1+j \right)}
= \,2\,\pi  \overset{\infty}{\underset{j =1}{\sum}}\; {\mathrm e}^{-b\,j^{2\,n +1}}
-2\,\pi \,b^{-{1}/({2\,n +1})}\,\Gamma \! \left(\frac{2\,n +2}{2\,n +1}\right)
\label{T2Gen}
\end{align}

Notice that if the index $j$ in the left-hand sum is odd, the next term in the series vanishes when $j$ is even and the left-hand side of \eqref{T2Gen} does not change. Therefore the integral is invariant when $N\rightarrow N+1$ if $N>0$ is odd.
That is
\begin{align} \nonumber
&\int_{-\infty}^{\infty}\zeta \! \left(\left(n +\frac{1}{2}\right) \left(2\,i\,v -4\,N +1\right)\right) b^{-i\,v}\,\Gamma \! \left(-2\,N +\frac{1}{2}+i\,v \right)d v \\
 = 
&b\int_{-\infty}^{\infty}\zeta \! \left(\left(n +\frac{1}{2}\right) \left(2\,i\,v -4\,N -1\right)\right) b^{-i\,v }\,\Gamma \! \left(-2\,N -\frac{1}{2}+i\,v \right)d v
\label{Tinv}
\end{align}
There are two interesting cases here, the first corresponding to $n=0$ -- see subsection (\ref{ssec:Seqm1}) -- when $N\rightarrow\infty$ so that
\begin{align}
\underset{N \rightarrow \infty}{\mathrm{lim}}\;\int_{-\infty}^{\infty}&\zeta \! \left(-N -\frac{1}{2}+i\,v \right) b^{N +\frac{1}{2}-i\,v}\,\Gamma \! \left(-N -\frac{1}{2}+i\,v \right)d v \\&+2\,\pi  \overset{\infty}{\underset{j =0}{\sum}}\; \frac{\zeta \! \left(-j \right) \left(-b \right)^{j}}{\Gamma \! \left(1+j \right)} 
 = 
2\,\pi  \overset{\infty}{\underset{j =1}{\sum}}\ {\mathrm e}^{-b\,j}-\frac{2\,\pi}{b}
\label{T2eLIM}
\end{align}
But, from \cite[Eqs. (24.2.1) and (25.6.3)]{NIST} 
\begin{equation}
\overset{\infty}{\underset{j =0}{\sum}}\; \frac{\zeta \! \left(-j \right) \left(-b \right)^{j}}{\Gamma \! \left(1+j \right)}
 = 
\frac{1}{b}\overset{\infty}{\underset{j =1}{\sum}}\; \frac{B_{j}\,b^{j}}{\Gamma \left(1+j \right)}
 = \frac{1}{e^{b}-1}-\frac{1}{b}
\label{But}
\end{equation}
and the elementary relation
\begin{equation}
\overset{\infty}{\underset{j =1}{\sum}}\; {\mathrm e}^{-b\,j} = 
\frac{1}{{\mathrm e}^{b}-1}
\label{Sexp}
\end{equation}
we find
\begin{equation}
\underset{N \rightarrow \infty}{\mathrm{lim}}\;\int_{-\infty}^{\infty} \zeta \! \left(-N -\frac{1}{2}+i\,v \right) b^{N +\frac{1}{2}-i\,v}\,\Gamma \! \left(-N -\frac{1}{2}+i\,v \right)d v=0\,.
\label{T2eLimz}
\end{equation}

\section{Summary}

It has been shown by way of a limited number of examples that summing over the free variable introduced by the inverse Mellin transform yields a number of interesting identities, each of which can be studied and pursued on their own. Some of these are possibly new and at least one (i.e. \eqref{Jb}) generalizes two classic identities due to Riemann.\newline

Of particular ongoing interest is the fact that the modified inverse Mellin transform studied here yields a contour integral that can then be transformed in such a way as to generate infinite series and integrals that can in turn be modified to produce unexpected identities involving hyperpowers. It is suggested that further study along the lines presented here is warranted. A cursory scan of \cite[Table 6]{IxForms} finds a plethora of Mellin transform pairs involving the fundamental functions of classical analysis and the most common hypergeometric functions (e.g. \cite{SQJ}), each of which possesses known transformations that could be invoked to generate new identities in the same manner as has been done here. For the ambitious reader, here is a suggestion for further emulation:


{\bf If $f(x)=\sin(x)$ and $F(s)=\sin(\pi\,s/2)\,\Gamma(s)$, from Theorem \ref{sec:T2} we obtain}
\begin{equation}
\overset{\infty}{\underset{k =1}{\sum}}\; \sin \! \left({k^{-s}}\right)
=\overset{\infty}{\underset{k =0}{\sum}}\; \frac{\left(-1\right)^{k}\,\zeta \! \left(\left(2\,k +1\right) s \right)}{\Gamma \! \left(2\,k +2\right)} \quad \Re(s)>1,
\label{Ex1}
\end{equation}
reducing to
\begin{equation}
\overset{\infty}{\underset{k =1}{\sum}}\; \left(\sin \! \left(\frac{1}{k}\right)-\frac{1}{k}\right)
=
\overset{\infty}{\underset{k =1}{\sum}}\; \frac{\left(-1\right)^{k}\,\zeta \! \left(2\,k +1\right)}{\Gamma \! \left(2\,k +2\right)}
\label{Ex1AB}
\end{equation}
if $s=1$. The possibilities are endless.

\section{Declarations}

\subsection{Compliance with Ethical Standards}
-Disclosure of potential conflicts of interest(none)

\subsection{Competing Interests - none}

The authors did not receive support from any organization for the submitted work. The second author was refused permission to apply for Open Access publication support from the Natural Sciences and Engineering Research Council of Canada.

\begin{appendices}
\appendix

\section{Proof of \eqref{eqno(6)}} \label{sec:Proof1}
\counterwithin*{equation}{section} 
\renewcommand\theequation{A.\arabic{equation}}

By evaluating the residues in \eqref{eqno(6)} as the contour is moved leftwards, we arrive at the following sum and its representation

\begin{equation}
\overset{\infty}{\underset{k =1}{\sum}}\; \frac{b^{2\,k -2}\,\zeta \! \left(1-2\,k \right)}{\Gamma \! \left(2\,k-1 \right)}
 = 
\,-\frac{1}{2}\overset{\infty}{\underset{k =1}{\sum}}\; \frac{b^{2\,k -2}\,B_{2\,k}}{ k\,\Gamma \! \left(2\,k -1\right)}\,,
\label{Bsum1}
\end{equation}
where \cite[Eq. 25.6.3]{NIST}
\begin{equation}
\zeta(-k)=\,(-1)^k\,B_{1+k}/(1+k)
\label{ZrefB}
\end{equation}
has been used, $B_{k}$ are Bernoulli numbers and we note that $\zeta(-2k)=0$. Following the application of \cite[Eq. 24.7.4]{NIST}
\begin{equation}
B_{2\,k} = 
\left(-1\right)^{k}\,\pi  \int_{0}^{\infty}\frac{t^{2\,k}}{\sinh^{2} \! \left(\pi \,t \right)}d t, 
\label{B2k}
\end{equation}
we now invert the sum and integral (both convergent) and, since

\begin{equation}
\overset{\infty}{\underset{k =1}{\sum}}\; \frac{\left(-1\right)^{k} \left(b\,t \right)^{2\,k}}{k\,\Gamma \! \left(2\,k-1 \right) }
 = 
\,-2\,\sin \! \left(b\,t \right) b\,t -2\,\cos \! \left(b\,t \right)+2
\label{Sumx}
\end{equation}
courtesy of \cite[Maple]{Maple23}, we are now left with the following integrals

\begin{equation}
\int_{0}^{\infty}\frac{1-\cos \; \left(b\,t \right)}{\sinh \! \left(\pi \,t \right)^{2}}\,d t
 = \frac{b\,\coth \! \left(\frac{b}{2}\right)-2}{2\,\pi}
\label{J3a}
\end{equation}
 and
 \begin{equation}
 \int_{0}^{\infty}\frac{\sin \! \left(b\,t \right) t}{\sinh \! \left(\pi \,t \right)^{2}}d t
 = 
\frac{b -\sinh \! \left(b \right)}{2\,\pi  \left(1-\cosh \! \left(b \right)\right)}
\label{J3b}
 \end{equation}
both courtesy of \cite[Mathematica]{Math23}. Putting it all together yields \eqref{eqno(6)}.

\section {Proof of \eqref{eqno(10)}} \label{sec:Proof2}

\renewcommand\theequation{B.\arabic{equation}}
From \eqref{eqno(10)} we are interested in the summation
\begin{equation}
\overset{\infty}{\underset{k =0}{\sum}}\; \frac{\zeta \! \left(-2\,k -1\right)}{\Gamma \! \left(2+2\,k \right) \left(a^{2\,k +1}-1\right)}
 = \,-\overset{\infty}{\underset{k =0}{\sum}}\; \frac{B_{2+2\,k}}{\Gamma \! \left(3+2\,k \right) \left(-1+a^{2\,k +1}\right)}
\label{T4b}
\end{equation}
employing \eqref{ZrefB}. Noting the identity \eqref{B2k}
and the well-known expansion
\begin{equation}
\frac{1}{1-{a^{-2\,k -1}}} = 
\overset{\infty}{\underset{j =0}{\sum}}\; \left(\frac{1}{a^{2\,k +1}}\right)^{j}, ~~a>1
\label{Aik}
\end{equation}
eventually leads us to consider
\begin{equation}
\overset{\infty}{\underset{k =0}{\sum}}\; \frac{\left(-1\right)^{k} \left({t}/{a}\right)^{2+2\,k}}{a^{j \,\left(2\,k +1\right)}\,\Gamma \! \left(3+2\,k \right)}
 = a^{j} \left(1-\cos \! \left(\frac{t}{a^{j +1}}\right)\right)
\label{Sid}
\end{equation}
after interchanging the order of summation. Now applying the identity \eqref{J3a} will eventually yield \eqref{eqno(10)}.

\end{appendices}

\bibliographystyle{unsrt}

\bibliography{c://physics//biblio}


\end{document}